\documentclass[preprint]{imsart}

\usepackage{amsthm,amsmath}
\RequirePackage[colorlinks,citecolor=blue,urlcolor=blue]{hyperref}

\usepackage[square,numbers]{natbib}

\startlocaldefs
\endlocaldefs

\begin{document}

\begin{frontmatter}

\title{On L\'evy's Brownian Motion and White Noise Space on the Circle}
\runtitle{L\'evy's Brownian Motion and White Noise Space on the Circle}

\begin{aug}

\author{\fnms{Chunfeng} \snm{Huang}\corref{}\ead[label=e1]{huang48@indiana.edu}},
\author{\fnms{Ao} \snm{Li} \ead[label=e2]{liao@indiana.edu}}

\address{Department of Statistics, Indiana University\\ Bloomington, IN 47408, USA }
 \affiliation{Indiana University}

\printead{e1,e2}
\phantom{E-mail:\ }



\end{aug}

\begin{abstract}
In this article, we show that the Brownian motion on the circle constructed by L\'evy (1959) is a regular Euclidean Brownian motion on the half-circle with its own mirror image on the other half-circle, and is degenerated in the sense of Minlos (1959). This raises the question of what the white noise is on the circle. We then formally define the white noise space and its associated Brownian bridge on the circle. 
\end{abstract}



\end{frontmatter}


\section{Introduction.}  The study of Brownian motions leads to many interesting and important branches in probability (\citet{Wiener1930}, \citet{Levy1948}, \citet{Gross1965}, \citet{Hida1970}). Most of these developments are in Euclidean spaces. On the circle, 
L\'evy (1959)\citep{Levy1959} constructed a Brownian motion $B(t)$ on the unit circle $S$ with the covariance function
\begin{equation} \label{eq:eq2}
\mbox{cov}(B(t), B(s)) = \frac{1}{2} (r(o, t) + r(o,s) - r(s,t)), \quad s ,t \in S,
\end{equation}
where $r(s,t)$ is the circular distance between $s$ and $t$, $o$ is an arbitrary pre-chosen point on $S$,  and $t \in [0,1]$ with $S = R/\mathcal{Z}$. In this construction, L\'evy \citep{Levy1959} assumes the existence of a white noise process on the circle first, and obtains the process $B(t)$ at $t$ to be the integration of this white noise process along the half-circle centered at $t$. Direct computation can show that this resulting process possesses the covariance function of (\ref{eq:eq2}). Additionally, it can be shown that such process has an equality (\citep{Levy1959}, equation 8):
\begin{equation} \label{eq:eq3}
B(t) + B(t') = B(s) + B(s'), \quad a.s., 
\end{equation}
where $t', s'$ are the opposites of $t,s$ on the circumference. Based on L\'evy's construction, it is obvious that $B(t)+B(t')$ is the integration of the white noise process along the entire circumference. Note that equation (\ref{eq:eq3}) implies that if this Brownian motion is given on the half-circle $\{ B(t), t \in [0, \frac{1}{2}] \}$, then for any point $s \in (\frac{1}{2}, 1)$, $B(s)$ can be obtained directly by equation (\ref{eq:eq3}):
\[
B(s) = B(0) + B(1/2) - B(s-1/2), \quad a.s..
\]
Note that for $s, t \in [0, 1/2]$, this Brownian motion has the covariance function of $\mbox{cov} (B(t), B(s)) = \min(s,t)$. In other words, this Brownian motion is actually a regular Euclidean Brownian motion on $[0,1/2]$, and the direct mirror image on $(1/2, 1]$. This poses a fundamental questions of the formal treatment of the white noise and Brownian type processes on the circle.   \\

Note that the white noise and other random processes on the circle have wide applications in probability and statistics, for example, the positive definite functions (\citet{Wood1995}, \citet{Pinkus2004}, \citet{Menegattoetal2006}, \citet{Sun2004}, \citet{Huangetal2011}, \citet{Gneiting2013}), the study of the intrinsic random functions (\citet{Huangetal2016}), the shape modeling (\citet{Hobolthetal2003, AlettiandRuffini2017}), and stochastic differential equations  (\citet{Lindgrenetal2011}). All these research are deeply rooted in the understanding of the white noise. The white noise is a mathematical idealization, and is not an ordinary random process (\citet{Yaglom1987}). It belongs to the generalized random process (\citet{Ito1953, GelfandandVilenkin1964}). Kuo (1996)\citep{kuo1996} offers several rigorous views of the white noise. One popular interpretation fo the white noise is to view it as a generalized derivative of the Brownian motion, and the existence and the construction of the Brownian motion in Euclidean space are well-developed (\citep{Wiener1930, Levy1948, Gross1965, Hida1970}).  However, as we just discussed, L\'evy's construction of Brownian motion on the circle does not appear to be appropriate for such understanding.  \\

In this manuscript, we construct the white noise space on the circle following the standard Gel'fand's triple and Minlos-Bochner theorem in Section 2. The Sobolev space of negative index ($-1$ in this article) for periodic functions $H_{-1}(S)$ is used in this construction. It is known that this space is the dual of the Sobolev space of order $1$, $H_{1,0}(S)$, with an additional condition of $\int_S f(t) dt = 0$. As we discuss in Section 3, this additional constraint leads naturally to the Brownian bridge, instead of Brownian motion on the circle.  \\

While the construction of the white noise space can follow the standard Gel'fand triple and Minlos-Bochner theorem directly and may seem to be rudimentary, they are presented in part of this manuscript for the completion. It helps understand the Brownian motion constructed in L\'evy is problematic, and reveals that the Brownian bridge is the natural choice in this situation.  \\

\section{White noise space on the circle}
The study of the white noise space can be found in  \citet{Hida1970}, \citet{Hidaetal1993}, \citet{kuo1996} and \citet{Si2012} based on the Bochner-Minlos Theorem. Minlos (1959, \citep{Minlos1959}) extended Bochner Theorem to infinite dimensional space with \citet{Kolmogorov1959}'s refinement in proof. A comprehensive discussion can be found in \citet{Hida1970}. Most of these developments are in Euclidean spaces. In this Section, we adapt this approach and formally define the white noise on the circle. 

We start with a Gel'fand triple for functions on the unit circle $S$:
\[
E \subset L^2(S) \subset E',
\]
where the space $E$ and its dual space $E'$ will be introduced shortly, and $L^2(S)$ is the collection of functions $f(t), t \in S$ such that
\[
\int_S | f(t) |^2 d t < \infty.
\]
This $L^2(S)$ is a Hilbert space with the inner product
\[
(f, g) = \int_S f(t) \overline{g}(t) d t,
\]
the associated norm
\[
\| f \| = \left( \int_S | f(t) |^2 d t \right)^{1/2},
\]
and the orthonormal basis
\[
\{ \psi_n (t) = e^{ i 2 \pi n t}, t \in S \}, \quad n \in \mathcal{Z}.
\]
For the space $E$, we consider the Sobolev space of order $1$ on the circle
\[
H_1(S) = \{ f \in L^2(S): \sum_n | (f, \psi_n)|^2 n^2 < \infty \}.
\]
Clearly $H_1(S) \subset L^2(S)$. Note that the Sobolev space of negative index is introduced as the dual space of the Sobolev space of the corresponding positive index with additional constraints (See \citet[Chap. 3]{Adams1975}). For Sobolev spaces of negative index for periodic functions, \citet[Chap. 5]{Robinson2001} defines $H_{-1}(S)$ as the dual space of $H_{1,0}(S)$, where
\begin{equation} \label{eq:H10}
H_{1,0}(S) = \{ f \in H_1(S): (f, \psi_0) = 0 \},
\end{equation}
then,
\[
H_{-1}(S) = (H_{1,0}(S))'.
\]
Note that $\psi_0(t) \equiv 1$, then, this extra condition is:
\[
(f, \psi_0) = \int_S f(t) d t,
\]
that is, the functions in $H_{1,0}(S)$ has the Fourier coefficient for $\psi_0$ to be zero. \\

{\bf Remark 1.} In Sobolev space literature, it is well-known that the space with negative index is defined as the dual with an additional constraint (\citep{Adams1975,Robinson2001}). As shown in Section 3, the constraint here with $(f, \psi_0) = 0$ plays an important role in obtaining Brownian bridge on the circle instead of Brownian motion that L\'evy constructed \citep{Levy1959}. \\

We now have the Gel'fand triple
\[
H_{1,0}(S) \subset L^2(S) \subset H_{-1}(S),
\]
with a natural pair
\[
(\xi, x), \quad \xi \in H_{1,0}(S), \quad x \in H_{-1}(S).
\]

{\bf Lemma 1.} The injection $H_{1,0}(S) \mapsto L^2(S)$ is of Hilbert-Schmidt type, so is the injection $L^2(S) \mapsto H_{-1}(S)$.

To show this, note that the norms in both $H_{1,0}(S)$ and $L^2(S)$ are
\[
\| f \|^2_{H_{1,0}(S)} = \sum_n n^2 | (f,\psi_n)|^2, \quad \| f \|^2 = \sum | f(, \psi_n)|^2.
\]
Let $ \xi_n = \frac{1}{n} \psi_n, n \ne 0$, then $\{ \xi_n \}$ is an orthonormal basis for $H_{1,0}(S)$, and
\[
\sum \| \xi_n \|^2 = \sum \frac{1}{n^2} < \infty.
\]
So the injection $H_{1,0}(S) \mapsto L^2(S)$ is of Hilbert-Schmidt type. Similarly, the injection $L^2(S) \mapsto H_{-1}(S)$ is also of Hilbert-Schmidt type. Parallel results can be found in \citep{kuo1996, Si2012}. \\

Consider the following functional $C(\xi), \xi \in H_{1,0}(S)$, 
\begin{equation} \label{eq:1}
C(\xi) = \exp \left( -\frac{1}{2} \| \xi \|^2 \right), \quad \xi \in H_{1,0}(S),
\end{equation}
where $\| \cdot \|$ is the $L^2(S)$ norm $\| \xi \|^2 = (\xi, \xi) = \int_S |\xi(t) |^2 d t$.\\

{\bf Lemma 2.} $C(\xi)$ is positive definite for $\xi \in H_{1,0}(S)$.

To show that $C(\xi)$ is positive definite, i.e., for any $z_j \in \mathcal{C}$ and $\xi_j \in H_{1,0}(S),  j=1,\ldots, n$,
\[
\sum_{j,k=1}^n z_j C(\xi_j - \xi_k) \overline{z}_k \ge 0.
\]
For $\xi_1, \ldots, \xi_n \in H_{1,0}(S)$, let $V$ be the (finite-dimensional) subspace of $H_{1,0}(S)$ spanned by $\xi_1, \ldots, \xi_n$ with the $L^2(R)$ norm. Let $\mu_V$ be the standard Gaussian measure on $V$, since $V$ is finite-dimensional, one can directly apply Bochner's Theorem. Then for any $\xi \in V$,
\[
\int_V e^{ i (\xi, y)} d \mu_V(y) = e^{ -\frac{1}{2} \| \xi \|^2}.
\]
Therefore,
\begin{eqnarray*}
&\,& \sum_{j,k=1}^n z_j C(\xi_j - \xi_k) \overline{z}_k = \sum_{j,k=1}^n \int_V z_j e^{ i (\xi_j - \xi_k, y)} d \mu_V(y) \overline{z}_k \\
&\,& = \int_V \left| \sum_{j=1}^n z_j e^{ i (\xi_j, y)} \right|^2 d \mu_V(y) \ge 0.
\end{eqnarray*}
The function $C(\xi)$ is obviously continuous on $H_{1,0}(S)$ and $C(0) = 1$. That is, $C(\xi)$ is a characteristic functional of $\xi \in H_{1,0}(S)$. \\

It is obvious that $C(\xi)$ is continuous on $H_{1,0}(S)$ and $C(0) = 1$. That is, $C(\xi)$ is a characteristic functional of $\xi \in H_{1,0}(S)$. Therefore, by Minlos' Theorem (\citet{Minlos1959, Kolmogorov1959}), we have the following theorem.

{\bf Theorem 3.} The functional $C(\xi), \xi \in H_{1,0}(S)$ in (\ref{eq:1}) determines the probability measure $\mu$ on $(H_{-1}(S), \mathcal{B})$ satisfying
\[
\int_{H_{-1}(S)} e^{ i (x, \xi)} d \mu(x) = C(\xi),
\]
where $\mathcal{B}$ is the Borel $\sigma$-algebra on $H_{-1}(S)$, i.e., the $\sigma$-algebra generated by the cylinder subsets of $H_{-1}(S)$. $C(\xi)$ is the characteristic functional of the measure $\mu$. \\

{\bf Definition 4.} The measure space $(H_{-1}(S), \mathcal{B}, \mu)$ determined by the characteristic functional $C(\xi)$ in (\ref{eq:1}) is called a while noise space on the circle. The measure $\mu$ is called the standard Gaussian measure on $H_{1,0}(S)$. \\

Notation wise, we denote the natural dual space pair $(\xi, x), \xi \in H_{1,0}(S), x \in H_{-1}(S)$ as $\xi(x)$, which is a random variable in this white noise space and has the following probability properties.

{\bf Proposition 5.} The random variable $\xi(x) \sim N(0, \| \xi \|^2)$, and the random variables $\xi_1(x)$ and $\xi_2(x)$ are independent when $(\xi_1, \xi_2) = 0$.

To show this, we first show that $(\xi_1(x), \ldots, \xi_n(x))$ follow a multivariate normal distribution through the following characteristic function derivation. 
\begin{eqnarray*}
&\,& \mbox{E} \exp \{ i t _1 \xi_1(x) + i t_2 \xi_2 (x) + \ldots + i t_n \xi_n (x) \} \\
&\,& = \int_{H_{-1}} \exp \{ i (t_1 \xi_1 + t_2 \xi_2 + \ldots + t_n \xi_n, x) \} d \mu(x) \\
&\,& = \exp \left( - \frac{1}{2} \| t_1 \xi_1 + t_2 \xi_2 + \ldots + t_n \xi_n\|^2 \right). 
\end{eqnarray*}
Therefore, $(\xi_1(x), \ldots, \xi_n(x))$ forms a Gaussian system. In particular, 
\[
\xi(x) \sim N(0, \| \xi \|^2),
\]
that is, $\xi(x)$ is a Gaussian random variable with mean zero, and variance $\| \xi \|^2$. And $\xi_1(x), \xi_2(x)$ are jointly bivariate Gaussian random vector with characteristic function
\[
\mbox{E} \exp \{ i t_1 \xi_1(x) + i t_2 \xi_2(x) \} = \exp \left( - \frac{1}{2} ( t_1^2 \| \xi_1 \|^2 + t_2^2 \| \xi_2 \|^2 + 2 t_1 t_2 (\xi_1, \xi_2) ) \right).
\]
When $(\xi_1, \xi_2) = 0$, $\xi_1(x)$ and $\xi_2(x)$ are independent.  \\

This Proposition indicates that 
\[
C(\xi_1 + \xi_2) = C(\xi_1) C(\xi_2),
\]
whenever
\[
\mbox{supp} (\xi_1) \cap \mbox{supp} (\xi_2) = \emptyset,
\]
where $\mbox{supp}(\xi)$ stands for the support of $\xi$. Based on the discussion in \citet[Sect. 4.5]{Hida1970}, this process has independent values at every moment, or equivalently, $\xi_1(x)$ and $\xi_2(x)$ are independent random variables if $\mbox{supp}(\xi_1) \cap \mbox{supp}(\xi_2) = \emptyset$. This is the so called white noise property, see also in  \citet{AdlerandTaylor2007, Lindgrenetal2011}. \\

\section{Brownian bridge on the circle.} In this Section, we formally introduce Brownian bridge and then further discuss the issues with L\'evy's Brownian motion construction in \citep{Levy1959}. 

Given the white noise space on the circle $(H_{-1}(S), \mathcal{S}, \mu)$ with the characteristic functional $C(\xi)$ in ({\ref{eq:1}). We have the following extension theorem.

{\bf Lemma 6. } The process $\xi(x), \xi \in H_{1,0}(S)$ has an extension to $L_0^2(S)$, where
\[
L_0^2(S) = \{ f \in L^2(S), (f, \psi_0) = 0 \}, 
\]
additionally, 
\[
\int_{H_{-1}(S)} (f, x) \overline{(g, x)} d \mu(x) = (f, g) 
\]

With $H_{1,0}(S)$ dense in $L_0^2(S)$, suppose $f \in L_0^2(S)$, take a sequence $\{ \xi_n \}$ in $H_{1,0}(S)$ such that $\xi_n \to f$ in $L_0^2(S)$. Then, the sequence $\{ \xi_n(x) \}$ of random variables is Cauchy in $L^2(H_{-1}(S), \mathcal{B}, \mu)$ since
\[
\int_{H_{-1}(S)} | \xi_k(x) - \xi_n(x) |^2 d \mu(x) = \| \xi_k - \xi_n \|^2.
\]
Denote
\[
f(x) = \lim_{n \to \infty} \xi_n(x), \mbox{ in } L^2(H_{-1}(S), \mathcal{B}, \mu).
\]
This limit is independent of the sequence. Moreover, the random variable $f(x)$ is normally distributed with mean $0$ and variance $\| f \|^2$, and for two $f, g \in L_0^2(S)$, the covariance of $f(x)$ and $g(x)$ is their $L_0^2(S)$ inner product by the continuity of extension, namely,
\[
\int_{H_{-1}(S)} f(x) g(x) d \mu(x) = (f, g),
\]
and Lemma 6 is proved. \\

In Euclidean space $R$, a Brownian motion can be introduced \citep{kuo1996} through $1_{[0,t)}$   as  
\[
(1_{[0,t)}, 1_{[0,s)}) = \min(s,t),
\]
which directly leads to a Brownian motion in $R$. However, $1_{[0,t)} \notin L_0^2(S)$ as it has Fourier coefficient $(1_{[0,t)}, \psi_0) = t \ne 0$, and cannot be extended from $H_{1,0}(S)$. 
Therefore, instead of $1_{[0,t)}$, we have
\[
1_{[0,t)} - t \in L_0^2(S),
\]
which yields
\[
(1_{[0,t)} - t, 1_{[0,s)} - s) = \min(s,t) - st.
\]
This is the covariance function for the Brownian bridge. Therefore, we have the following result.

{\bf Theorem 7.} Given the white noise space $(H_{-1}(S), \mathcal{B}, \mu)$ on the circle, the generalized random process
\begin{equation} \label{eq:4}
(1_{[0,t)} -t, x), \quad x \in H_{-1}(S),
\end{equation}
is a Brownian bridge. \\

{\bf Remark 2}. On real line, $1_{[0,t)}$ is usually used to obtain Brownian motion \citep{kuo1996}. It is clear that the additional condition of $(f, \psi_0)=1$ of the test function space $H_{1,0}(S)$ introduced in Section 2 plays a direct role in obtaining this Brownian bridge. 
The constraint of $(f, \psi_0)=0$ in equation (\ref{eq:H10}) clearly indicate that Brownian bridge is a natural result on the circle instead of Brownian motion, see also \citep{Hidaetal1993}.  \\

{\bf Remark 3.} As discussed in the Introduction Section, L\'evy's Brownian motion on the circle (\citep{Levy1959}, \citep{Gangolli1967}) is a regular one-dimensional Euclidean Brownian motion on $[0,1/2]$, and mere mirror imagine on $(1/2, 1]$.  Given our discussion of the white noise space on the circle, one can obtain that there is a test function $\eta_t(u), u \in S$ with
\[
\eta_t(u) = \frac{1}{\sqrt{2} \pi} \sum_{k \mbox{ is odd}} \frac{ e^{i 2 \pi k t}-1}{|k|} e^{i 2 \pi k u}.
\]
It is clear that $\eta_t(u) \in L^2(S)$ and $(\eta_t, \psi_0) = 0$. Lemma 6 can apply and direct computation shows that
\[
(\eta_t, \eta_s) = \min(s,t).
\]
That is, the generalized random process
\[
\eta_t(x), \quad x \in H_{-1}(S)
\]
shares the same covariance with the Brownian motion in \citet{Levy1959}. It is also obvious from this development that the function $\eta_t(u)$ only has odd terms, and the process becomes degenerated in the sense of \citet[Section 3]{Minlos1959}. The covariance in (\ref{eq:eq2}) can be shown to be not strictly positive definite on the circle (\citet{Pinkus2004}). \\

{\bf Remark 4.}  In Euclidean spaces, a white noise, sometimes, can be thought of as the generalized derivative of a Brownian motion. With our discussion, L\'evy's Brownian motion on the circle is not proper for this consideration. However, the Brownian bridge in (\ref{eq:4}) is well-defined and its generalized derivative can be viewed as the white noise on the circle.


\bibliography{mybibfile}





\end{document}